# The Parametric Instability of Motion at Resonance as a Source of Chaotic Behaviour at the Restricted Three body Problem

A.E. Rosaev

*FGUP NPC Nedra, Yaroslavl, Russia*

The planetary restricted three-body problem (RTBP) is considered. The primary mass M is much more than another mass *m*, which revolves around M. The massless probe particle $m_0$ moves on elliptic orbit, is perturbed by *m*.

Our target is - to study stability and instability in motion of *m* at resonance. The resonant type of motion is widespread in the Solar system. However, investigation of resonance, in particular, calculation of width of unstable regions, is a difficult problem. It is well known, that equations of motion of restricted 3-body problem in rotating rectangular frame may be reduced to the second order differential equation with periodic coefficients (Hill's equation). Here we derive Hill's equation at cylindrical coordinate frame. It gives the ability to estimate width and position of the unstable zones. Moreover, the dependence of the position of instable zones from orbital eccentricity of the probe particle is derived. This result may be compared with distribution of asteroids near the resonance. It is noted, that overlapping of instability zones in strongly perturbed system can be a source of the chaotic behaviour.

Initially, we use some simplifications - planar motion and circular orbits of the perturbing mass $m_i$.

## INTRODUCTION. MAIN EQUATIONS

There are some possible coordinate frames, useful for perturbed motion consideration. On our view, cylindrical system has some advantages over other ones, because the variation of one of the coordinates – central distance R - is restricted forever and may be considered as a small parameter at the problem. The main equations for the planar circular 3-body problem in absolute frame are:

$$\frac{d^2 R(t)}{dt^2} - R(t)\left(\frac{d\lambda(t)}{dt}\right)^2 = \frac{dU}{dR}$$

$$\frac{d(R(t)^2 \lambda'(t))}{dt} = \frac{dU}{d\lambda}$$
(1)

where perturbation function in absolute frame:

$$U = \frac{Gm}{(r^2 + R^2 - 2rR\cos(S))^{1/2}} + \frac{GM}{R}$$
(2)

where M- mass of the primary; m - mass of the perturbing body; R, r - distance from the mass center test particle and perturbing body accordingly; S - angle between perturbed and perturbing body, G- constant of gravity, $\delta\lambda(t)=(\omega - \omega_s)t+ \varphi \equiv S$ - longitude of perturbed body, $\omega, \omega_s$ - mean motions perturbed and perturbing body accordingly, $\varphi$ - initial phase, and:

$$\Delta = \sqrt{R^2 + r^2 - 2Rr\cos(S)} \qquad (3)$$

The equations of planar restricted Hill's problem in rectangular frame:

$$\frac{d^2X}{dt^2} - 2m\frac{dY}{dt} + fX = 0$$
$$\frac{d^2Y}{dt^2} + 2m\frac{dX}{dt} + gY = 0 \qquad (4)$$

where $m$ – perturbing body mass, $X, Y$ – rectangle coordinates, $f$ and $g$ – known functions, $t$- time.

It may be reduced to a Hill equation for normal distance from variation orbit $x$ [1]:

$$\frac{d^2x}{dt^2} + \omega^2(t)x = f(t) \qquad (5)$$

where $f(t)$ – known function of time, $\omega(t)$ - periodic function of time (to be determined).

The cylindrical coordinate frame has some advantages for many applications of planetary dynamics due to a taking into account problem symmetry. The according equations in cylindrical coordinates may be derived from (5). However, they may be written independently. There are two small parameters in problem: x - radial shift from intermediate (in variations) orbit and r/R – the ratio of mean distance of perturbing and perturbed body. Let $R(t)=R+x(t)$, where R is constant. Accordingly, there are two ways of linearization (two possible sequence of expansion). And in both cases, the angle between perturbed and perturbing body must be accordingly expanded:

## THE EQUATION IN VARIATIONS

The principle part of expansion perturbation function may be written with using Legendre's polynomials $P_i(\cos\delta\lambda)$. In case of outer perturbing body [2]:

$$U = -\frac{\gamma M}{R} - \sum_{p=2}^{\infty} \frac{\gamma m}{r} y^p P_p(\cos(\delta\lambda)) \ , \qquad y=R/r \qquad (6)$$

$$\frac{d^2U}{dx^2} = \frac{2\gamma M}{R^3} - \frac{\gamma m}{r^3} \sum_{p=2}^{\infty} p(p-1) y^{p-2} P_p(\cos(\delta\lambda)) \qquad (7)$$

There are two ways of Hill (Mathieu) equation deriving. At linear approach by x we have:

$$\ddot{x} + \omega^2 x \approx f(t) \tag{8}$$

$$f(t) = \frac{L^2_0}{R_0^3} - \frac{\gamma M}{R_0^2} - \sum_{p=2}^{\infty} p \frac{\gamma m}{r^2} y_0^{p-1} P_p(\cos(\delta\lambda))$$

where ω depend from time:

$$\omega^2 = \omega_0^2 \left(1 + \sum h_n \cos n \varpi t\right) \tag{9}$$

where:

$$\omega^2 = \gamma M / R^3 - \gamma m / \Delta^3 (1 - 3(R - r\cos(\delta\omega t))/\Delta^2)) \tag{10}$$

Then the expansion by power y=$R_0$/r is applied. At zero order by y:

$$\omega^2 = \gamma M / R_o^3 + \gamma m / (2r^3) + 3/2 \gamma m / r^3 (\cos(2\delta\omega t)) \tag{11}$$

After using another sequence of expansion, first expansion by power y=$R_0$/r and then by power x:

$$\omega^2 = \frac{\gamma M}{R_0^3} + \frac{\gamma m}{r^3} \sum_{p=2}^{\infty} p(p-1) y^{p-2} P_p(\cos(\delta\lambda)) \tag{12}$$

where P(cos) is a Legendre polynomials. After the substitution:

$$P_2(\cos(S)) = 1/4(3(\cos(2S)+1)) \tag{13}$$

we can see that, at linear approach by x and zero order by y, ω=ω and both way of deriving ω are equivalent.

Now we can easy take into account the eccentricity of perturbed particle orbit and show, that such generalization leads to very interesting results.

At the elliptic osculating orbit, after the entering osculating orbital elements *a* - semimajor axis, *e* - eccentricity, approximately:

$$L_0 = \sqrt{\gamma M a(1-e^2)} \tag{14}$$

and $R_0$ depend from time. The phase-averaged value of $R_0$ approximately:

$$\overline{R}_0 = \frac{1}{2\pi}\int_0^{2\pi}\frac{a(1-e^2)}{1+e\cos(\varphi)}d\phi = a(1-e^2)^{1/2} \tag{15}$$

After angular moment substitution and linearization:

$$\frac{3L_0^2}{R_0^4} - \frac{2\gamma M}{R_0^3} \approx \frac{\gamma M}{a^3(1-e^2)^{3/2}}\left(1-3/2e^2\right) \approx \frac{\gamma M}{a^3}\left(1-\frac{3}{4}e^4\right) \tag{16}$$

It seems their dependence from eccentricity. So, if not *e=0*, mean motion must differ from keplerian circular. Finally, for the probe particle elliptic motion case:

$$\omega^2 = \frac{\gamma M}{a^3}\left(1-\frac{3}{4}e^4\right) + \frac{\gamma m}{r^3}\sum_{p=2}^{\infty}p(p-1)y^{p-2}P_p(\cos(\delta\lambda)) \tag{17}$$

Then we can group terms with equal *pλ* by using Laplace's coefficients. The main frequency may be performed:

$$\omega^2 = \frac{\gamma M}{R_0^3}\left(1-\frac{3}{4}e^4\right) + \frac{\gamma m}{r^3}\sum_{p=0}^{\infty}b_p(\cos(p\delta\lambda)) \tag{18}$$

where $b_p$ is easy to calculate numerically. For the periodic, or for the quasi-periodic solutions: $\delta\lambda = \sum \lambda_n \cos(n\varpi t)$. This is the Hill's equation:

$$\ddot{x} + \omega^2 x \approx f(t) \qquad \omega^2 = \omega_0^2\left(1+\sum h_n \cos n\varpi t\right) \tag{19}$$

where, at this performance, $\omega_0^2 = \frac{\gamma M}{R_0^3}\left(1-\frac{3}{4}e^4\right)$

Hill's equation also is linear (with floating factors), so its approximately solution possible to get a position of instability.

### THE STABILITY OF EQUATION IN VARIATIONS SOLUTION

For investigation of stability, *x* may be negligible small. From instability in this case, it will be followed instability in general. So, we can restrict only 1-st order in expansion by power *x*. The Mathieu's equation is a limit case of Hill's equation and it is more simply for studying. On the other side, because both equations (Hill and Mathieu) are linear, the main area of instability of Mathieu equation must be present in Hill's equation solution.

Consequently, and in this case exist the areas of instability, complying with condition:

$$\omega/\omega_s = nk/(nk - 2(1-\alpha/2)) \tag{20}$$

Thereby received that orbits near resonances (2n+1)/(2n-1) unstable parametric. This conclusion completely coincides with results of declared problem studies in paper Hadjidemetriou [3] by methods of matrix algebra.

Directly, from view of linear equations in elliptic case, we can explain one interesting feature - centres of resonant zones are shifted outside relative exact commensurability (resonance) in direction toward the source of perturbation. In simple case one perturbing body, centres of unstable resonant zones moved away from exact commensurabilities. So, exact commensurability may be outside relatively according unstable zone!

In result, we have a small drift of the centres of the instable zones from exact commensurabilities **toward** a source of perturbation. The width of unstable zones $\varepsilon$ for the first and second order is determined [4]:

$$\varepsilon_1 = \omega_0 \frac{nh_n}{2(n-2)} \qquad \varepsilon_{12} = \omega_0 \frac{nh^2_n}{8(n-2)} \qquad (21)$$

where values of $h_n$ proportional satellite mass (table 1):

$$h_p = \frac{\omega_s^2}{\omega_0^2} b_p = \frac{m/M}{(r/R_0)^3} \frac{1}{\left(1 - \frac{3}{4}e^4\right)} b_p \qquad (22)$$

Table 1

Width and position of different order resonance unstable zones

| order k | resonance | exact commensurability distance | b | Width r (e=0) |
|---|---|---|---|---|
| 3 | 3:1 | 2.501120 | 0.077800 | 0.0120 |
| 4 | 2:1 | 3.277395 | 0.107630 | 0.0114 |
| 6 | 3:2 | 3.700976 | 0.132076 | 0.0108 |

These expressions show, that, at strong perturbations, unstable resonance zones may be overlapped. The beginning from this, according with chaotic motion appearance, and when such zones fill all phase space, the behavior of system becomes completely chaotic.

The condition of overlapping is:

$$[(n+1)/(n-(1-3/4e^4)) - n/(n-(2-3/4e^4))] = \frac{nm/M b_n}{2(n-2)(1-3/8e^4)} \qquad (23)$$

It seems, that width of instability increased with eccentricity. On the other hand, distance between two close unstable areas ($\Delta\omega$) decreased. Width of areas of instability quickly decreases with the growing of order n. However, for each m, M, the value of n, since of which the overlapping of unstable zones take place, is present (Fig.3). For $e=0$:

$$n = \frac{mb \pm \sqrt{(mb)^2 + 16Mmb}}{2mb} \qquad (24)$$

and it is evident, that overlapping of unstable zones possible only at limit $n \rightarrow \infty$. It means, that effect of overlapping resonance instable zones and related chaotic behavior, strongly depends from the orbital eccentricity of probe particle $m$.

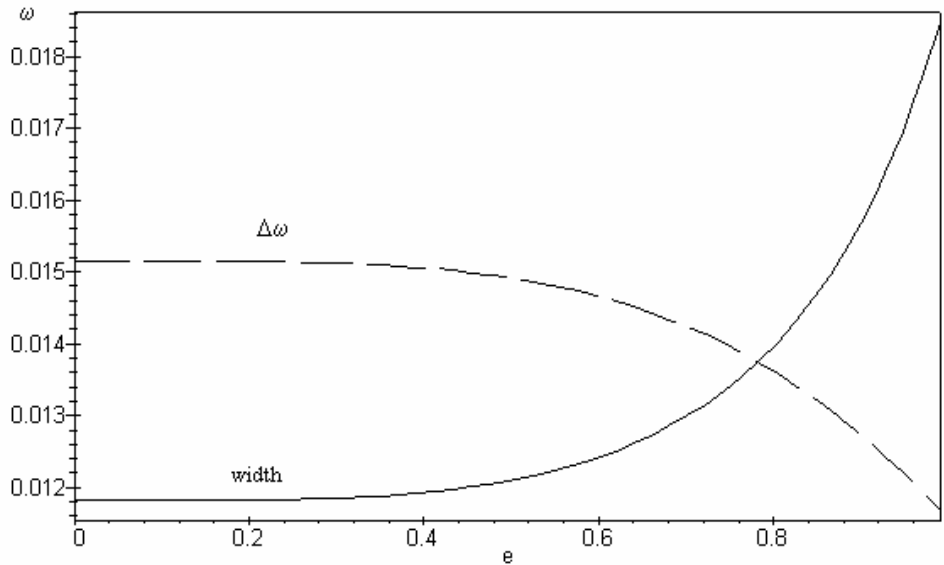

Fig.1. The dependence width and distance between unstable zones from eccentricity. m/M=0.1, n=13.

## RESULTS AND CONCLUSIONS

In result, we obtain very suitable equation in variation to study perturbed motion at resonance.

The main result of this work may be formulated in such form:

***Theorem 1.*** *There are not exist stable orbits with e=0 and i=0 in neighborhood of mean motion resonances (2n-1)/(2n+1).*

***Theorem 2.*** *Positions and width of unstable zones depends from eccentricity. It leads to overlapping of them at high order of resonance.*

Width of areas of instability quickly decreases with the growing of order n. As a result, to be unstable parametric first of all resonances 3:1, 2:1, 3:2, 5:3, in the case of circular orbits is not discovered instability near resonances 4:1, 5:1, 6:1, 5:2.

The lost of stability may be by two ways - by mean motion or by eccentricity change. It is seems, that second way is more probable, and leads to non-zero eccentricity. Evidently, that orbits with *e=0* in perturbed problem *no exist* not only in resonant case ("forced eccentricity"). It is means, that consideration elliptic motion of test particle is necessary.